\expandafter\chardef\csname pre amssym.def at\endcsname=\the\catcode`\@
\catcode`\@=11

\def\undefine#1{\let#1\undefined}
\def\newsymbol#1#2#3#4#5{\let\next@\relax
 \ifnum#2=\@ne\let\next@\msafam@\else
 \ifnum#2=\tw@\let\next@\msbfam@\fi\fi
 \mathchardef#1="#3\next@#4#5}
\def\mathhexbox@#1#2#3{\relax
 \ifmmode\mathpalette{}{\m@th\mathchar"#1#2#3}%
 \else\leavevmode\hbox{$\m@th\mathchar"#1#2#3$}\fi}
\def\hexnumber@#1{\ifcase#1 0\or 1\or 2\or 3\or 4\or 5\or 6\or 7\or 8\or
 9\or A\or B\or C\or D\or E\or F\fi}

\font\tenmsa=msam10
\font\sevenmsa=msam7
\font\fivemsa=msam5
\newfam\msafam
\textfont\msafam=\tenmsa
\scriptfont\msafam=\sevenmsa
\scriptscriptfont\msafam=\fivemsa
\edef\msafam@{\hexnumber@\msafam}
\mathchardef\dabar@"0\msafam@39
\def\dashrightarrow{\mathrel{\dabar@\dabar@\mathchar"0\msafam@4B}}
\def\dashleftarrow{\mathrel{\mathchar"0\msafam@4C\dabar@\dabar@}}

\def\ulcorner{\delimiter"4\msafam@70\msafam@70 }
\def\urcorner{\delimiter"5\msafam@71\msafam@71 }
\def\llcorner{\delimiter"4\msafam@78\msafam@78 }
\def\lrcorner{\delimiter"5\msafam@79\msafam@79 }
\def\yen{{\mathhexbox@\msafam@55 }}
\def\checkmark{{\mathhexbox@\msafam@58 }}
\def\circledR{{\mathhexbox@\msafam@72 }}
\def\maltese{{\mathhexbox@\msafam@7A }}

\font\tenmsb=msbm10
\font\sevenmsb=msbm7
\font\fivemsb=msbm5
\newfam\msbfam
\textfont\msbfam=\tenmsb
\scriptfont\msbfam=\sevenmsb
\scriptscriptfont\msbfam=\fivemsb
\edef\msbfam@{\hexnumber@\msbfam}

\catcode`\@=\csname pre amssym.def at\endcsname

\expandafter\ifx\csname pre amssym.tex at\endcsname\relax \else \endinput\fi
\expandafter\chardef\csname pre amssym.tex at\endcsname=\the\catcode`\@
\catcode`\@=11
\newsymbol\boxdot 1200
\newsymbol\boxplus 1201
\newsymbol\boxtimes 1202
\newsymbol\square 1003
\newsymbol\blacksquare 1004
\newsymbol\centerdot 1205
\newsymbol\lozenge 1006
\newsymbol\blacklozenge 1007
\newsymbol\circlearrowright 1308
\newsymbol\circlearrowleft 1309
\undefine\rightleftharpoons
\newsymbol\rightleftharpoons 130A
\newsymbol\leftrightharpoons 130B
\newsymbol\boxminus 120C
\newsymbol\Vdash 130D
\newsymbol\Vvdash 130E
\newsymbol\vDash 130F
\newsymbol\twoheadrightarrow 1310
\newsymbol\twoheadleftarrow 1311
\newsymbol\leftleftarrows 1312
\newsymbol\rightrightarrows 1313
\newsymbol\upuparrows 1314
\newsymbol\downdownarrows 1315
\newsymbol\upharpoonright 1316
 
\newsymbol\downharpoonright 1317
\newsymbol\upharpoonleft 1318
\newsymbol\downharpoonleft 1319
\newsymbol\rightarrowtail 131A
\newsymbol\leftarrowtail 131B
\newsymbol\leftrightarrows 131C
\newsymbol\rightleftarrows 131D
\newsymbol\Lsh 131E
\newsymbol\Rsh 131F
\newsymbol\rightsquigarrow 1320
\newsymbol\leftrightsquigarrow 1321
\newsymbol\looparrowleft 1322
\newsymbol\looparrowright 1323
\newsymbol\circeq 1324
\newsymbol\succsim 1325
\newsymbol\gtrsim 1326
\newsymbol\gtrapprox 1327
\newsymbol\multimap 1328
\newsymbol\therefore 1329
\newsymbol\because 132A
\newsymbol\doteqdot 132B
 
\newsymbol\triangleq 132C
\newsymbol\precsim 132D
\newsymbol\lesssim 132E
\newsymbol\lessapprox 132F
\newsymbol\eqslantless 1330
\newsymbol\eqslantgtr 1331
\newsymbol\curlyeqprec 1332
\newsymbol\curlyeqsucc 1333
\newsymbol\preccurlyeq 1334
\newsymbol\leqq 1335
\newsymbol\leqslant 1336
\newsymbol\lessgtr 1337
\newsymbol\backprime 1038
\newsymbol\risingdotseq 133A
\newsymbol\fallingdotseq 133B
\newsymbol\succcurlyeq 133C
\newsymbol\geqq 133D
\newsymbol\geqslant 133E
\newsymbol\gtrless 133F
\newsymbol\sqsubset 1340
\newsymbol\sqsupset 1341
\newsymbol\vartriangleright 1342
\newsymbol\vartriangleleft 1343
\newsymbol\trianglerighteq 1344
\newsymbol\trianglelefteq 1345
\newsymbol\bigstar 1046
\newsymbol\between 1347
\newsymbol\blacktriangledown 1048
\newsymbol\blacktriangleright 1349
\newsymbol\blacktriangleleft 134A
\newsymbol\vartriangle 134D
\newsymbol\blacktriangle 104E
\newsymbol\triangledown 104F
\newsymbol\eqcirc 1350
\newsymbol\lesseqgtr 1351
\newsymbol\gtreqless 1352
\newsymbol\lesseqqgtr 1353
\newsymbol\gtreqqless 1354
\newsymbol\Rrightarrow 1356
\newsymbol\Lleftarrow 1357
\newsymbol\veebar 1259
\newsymbol\barwedge 125A
\newsymbol\doublebarwedge 125B
\undefine\angle
\newsymbol\angle 105C
\newsymbol\measuredangle 105D
\newsymbol\sphericalangle 105E
\newsymbol\varpropto 135F
\newsymbol\smallsmile 1360
\newsymbol\smallfrown 1361
\newsymbol\Subset 1362
\newsymbol\Supset 1363
\newsymbol\Cup 1264
 
\newsymbol\Cap 1265
 
\newsymbol\curlywedge 1266
\newsymbol\curlyvee 1267
\newsymbol\leftthreetimes 1268
\newsymbol\rightthreetimes 1269
\newsymbol\subseteqq 136A
\newsymbol\supseteqq 136B
\newsymbol\bumpeq 136C
\newsymbol\Bumpeq 136D
\newsymbol\lll 136E
 
\newsymbol\ggg 136F
 
\newsymbol\circledS 1073
\newsymbol\pitchfork 1374
\newsymbol\dotplus 1275
\newsymbol\backsim 1376
\newsymbol\backsimeq 1377
\newsymbol\complement 107B
\newsymbol\intercal 127C
\newsymbol\circledcirc 127D
\newsymbol\circledast 127E
\newsymbol\circleddash 127F
\newsymbol\lvertneqq 2300
\newsymbol\gvertneqq 2301
\newsymbol\nleq 2302
\newsymbol\ngeq 2303
\newsymbol\nless 2304
\newsymbol\ngtr 2305
\newsymbol\nprec 2306
\newsymbol\nsucc 2307
\newsymbol\lneqq 2308
\newsymbol\gneqq 2309
\newsymbol\nleqslant 230A
\newsymbol\ngeqslant 230B
\newsymbol\lneq 230C
\newsymbol\gneq 230D
\newsymbol\npreceq 230E
\newsymbol\nsucceq 230F
\newsymbol\precnsim 2310
\newsymbol\succnsim 2311
\newsymbol\lnsim 2312
\newsymbol\gnsim 2313
\newsymbol\nleqq 2314
\newsymbol\ngeqq 2315
\newsymbol\precneqq 2316
\newsymbol\succneqq 2317
\newsymbol\precnapprox 2318
\newsymbol\succnapprox 2319
\newsymbol\lnapprox 231A
\newsymbol\gnapprox 231B
\newsymbol\nsim 231C
\newsymbol\ncong 231D
\newsymbol\diagup 231E
\newsymbol\diagdown 231F
\newsymbol\varsubsetneq 2320
\newsymbol\varsupsetneq 2321
\newsymbol\nsubseteqq 2322
\newsymbol\nsupseteqq 2323
\newsymbol\subsetneqq 2324
\newsymbol\supsetneqq 2325
\newsymbol\varsubsetneqq 2326
\newsymbol\varsupsetneqq 2327
\newsymbol\subsetneq 2328
\newsymbol\supsetneq 2329
\newsymbol\nsubseteq 232A
\newsymbol\nsupseteq 232B
\newsymbol\nparallel 232C
\newsymbol\nmid 232D
\newsymbol\nshortmid 232E
\newsymbol\nshortparallel 232F
\newsymbol\nvdash 2330
\newsymbol\nVdash 2331
\newsymbol\nvDash 2332
\newsymbol\nVDash 2333
\newsymbol\ntrianglerighteq 2334
\newsymbol\ntrianglelefteq 2335
\newsymbol\ntriangleleft 2336
\newsymbol\ntriangleright 2337
\newsymbol\nleftarrow 2338
\newsymbol\nrightarrow 2339
\newsymbol\nLeftarrow 233A
\newsymbol\nRightarrow 233B
\newsymbol\nLeftrightarrow 233C
\newsymbol\nleftrightarrow 233D
\newsymbol\divideontimes 223E
\newsymbol\varnothing 203F
\newsymbol\nexists 2040
\newsymbol\Finv 2060
\newsymbol\Game 2061
\newsymbol\mho 2066
\newsymbol\eth 2067
\newsymbol\eqsim 2368
\newsymbol\beth 2069
\newsymbol\gimel 206A
\newsymbol\daleth 206B
\newsymbol\lessdot 236C
\newsymbol\gtrdot 236D
\newsymbol\ltimes 226E
\newsymbol\rtimes 226F
\newsymbol\shortmid 2370
\newsymbol\shortparallel 2371
\newsymbol\smallsetminus 2272
\newsymbol\thicksim 2373
\newsymbol\thickapprox 2374
\newsymbol\approxeq 2375
\newsymbol\succapprox 2376
\newsymbol\precapprox 2377
\newsymbol\curvearrowleft 2378
\newsymbol\curvearrowright 2379
\newsymbol\digamma 207A
\newsymbol\varkappa 207B
\newsymbol\Bbbk 207C
\newsymbol\hslash 207D
\undefine\hbar
\newsymbol\hbar 207E
\newsymbol\backepsilon 237F
\catcode`\@=\csname pre amssym.tex at\endcsname

\magnification=1200
\hsize=468truept
\vsize=646truept
\voffset=-10pt
\parskip=1pc
\baselineskip=14truept
\count0=1

\dimen100=\hsize

\def\leftill#1#2#3#4{
\medskip
\line{$
\vcenter{
\hsize = #1truept \hrule\hbox{\vrule\hbox to  \hsize{\hss \vbox{\vskip#2truept
\hbox{{\copy100 \the\count105}: #3}\vskip2truept}\hss }
\vrule}\hrule}
\dimen110=\dimen100
\advance\dimen110 by -36truept
\advance\dimen110 by -#1truept
\hss \vcenter{\hsize = \dimen110
\medskip
\noindent { #4\par\medskip}}$}
\advance\count105 by 1
}
\def\rightill#1#2#3#4{
\medskip
\line{
\dimen110=\dimen100
\advance\dimen110 by -36truept
\advance\dimen110 by -#1truept
$\vcenter{\hsize = \dimen110
\medskip
\noindent { #4\par\medskip}}
\hss \vcenter{
\hsize = #1truept \hrule\hbox{\vrule\hbox to  \hsize{\hss \vbox{\vskip#2truept
\hbox{{\copy100 \the\count105}: #3}\vskip2truept}\hss }
\vrule}\hrule}
$}
\advance\count105 by 1
}
\def\midill#1#2#3{\medskip
\line{$\hss
\vcenter{
\hsize = #1truept \hrule\hbox{\vrule\hbox to  \hsize{\hss \vbox{\vskip#2truept
\hbox{{\copy100 \the\count105}: #3}\vskip2truept}\hss }
\vrule}\hrule}
\dimen110=\dimen100
\advance\dimen110 by -36truept
\advance\dimen110 by -#1truept
\hss $}
\advance\count105 by 1
}
\def\insectnum{\copy110\the\count120
\advance\count120 by 1
}

\font\ninerm=cmr9
\font\eightrm=cmr8

\font\tenrm=cmr10 at 10pt

\font\sc=cmcsc10

\def\msb{\fam\msbfam\tenmsb}

\def\bbc{{\msb C}}

\def\bbp{{\msb P}}
\def\bbq{{\msb Q}}

\def\bbz{{\msb Z}}
\def\grD{\Delta}

\def\grL{\Lambda}

\def\grl{\lambda}

\def\grz{\zeta}
\font\svtnrm=cmr17

\font\aa=eufm10

\def\Got#1{\hbox{\aa#1}}

\def\gsp1{{\Got s}{\Got p}(1)}

\def\mh-1{\hat{\mu}^{-1}(0)}
\def\n-1c{\nu^{-1}(c)}
\def\m-1{\mu^{-1}(0)}
\def\p-1{\pi^{-1}}
\def\p'-1{\prime{\pi}^{-1}}

\def\cald{{\cal D}}

\def\calf{{\cal F}}

\def\calm{{\cal M}}

\def\cals{{\cal S}}

\def\calz{{\cal Z}}

\def\la#1{\hbox to #1pc{\leftarrowfill}}
\def\ra#1{\hbox to #1pc{\rightarrowfill}}

\def\fract#1#2{\raise4pt\hbox{$ #1 \atop #2 $}}
\def\decdnar#1{\phantom{\hbox{$\scriptstyle{#1}$}}
\left\downarrow\vbox{\vskip15pt\hbox{$\scriptstyle{#1}$}}\right.}

\def\bowtie{\hbox to 1pt{\hss}\raise.66pt\hbox{$\scriptstyle{>}$}
\kern-4.9pt\triangleleft}
\def\hsmash{\triangleright\kern-4.4pt\raise.66pt\hbox{$\scriptstyle{<}$}}
\def\boxit#1{\vbox{\hrule\hbox{\vrule\kern3pt
\vbox{\kern3pt#1\kern3pt}\kern3pt\vrule}\hrule}}

\def\za{\vrule height6pt width4pt depth1pt}

\font\aa=eufm10

\def\Got#1{\hbox{\aa#1}}

\def\bfw{{\bf w}}

\def\cald{{\cal D}}

\def\calf{{\cal F}}

\def\calm{{\cal M}}

\def\cals{{\cal S}}

\def\calz{{\cal Z}}

\def\gG{{\Got G}}

\font\svtnrm=cmr17

\font\bsc=cmcsc10 at 10truept

\def\Se{Sasakian-Einstein }

\settabs 9\columns

\font\svtnrm=cmr17

\font\bsc=cmcsc10 at 10truept

\def\Se{Sasakian-Einstein }

\def\Div{\hbox{Div~}}
\def\r5int{0}

\centerline{\svtnrm New Einstein Metrics on $8\#(S^2\times S^3)$}

\bigskip
\centerline{\sc Charles P. Boyer~~ Krzysztof Galicki}
\footnote{}{\ninerm During the preparation of this work the authors 
were partially supported by NSF grants DMS-9970904 and DMS-0203219. 2000 
}
\bigskip

\centerline{\vbox{\hsize = 5.85truein
\baselineskip = 12.5truept
\eightrm
\noindent {\bsc Abstract:}
We show that $\scriptstyle{\#8(S^2\times S^3)}$ admits two 8-dimensional
complex families of inequivalent non-regular \Se structures. These are the first
known non-regular Sasakian-Einstein metrics on this 5-manifold.}} \tenrm

\bigskip
\baselineskip = 10 truept
\centerline{\bf  Introduction}  
\bigskip

Recently, the authors [BG2] 
introduced a new method for showing the existence of 
Sasakian-Einstein metrics on compact 
simply connected 5-dimensional spin manifolds. This 
method was based on work of Demailly and Koll\'ar [DK] who gave sufficient algebraic 
conditions on log del Pezzo surfaces anticanonically embedded in weighted projective spaces 
to guarantee the existence of a K\"ahler-Einstein orbifold metric. This, in turn enabled us to 
construct \Se metrics on certain $S^1$ V-bundles over these log del Pezzo surfaces. One 
could then use known monodromy techniques on the links of isolated hypersurface 
singularities together with a classification result of Smale [Sm] to identify the 5-manifold. The 
Demailly and Koll\'ar methods were further developed by Johnson and Koll\'ar [JK] where a 
computer code was written to solve the algebraic equations. The authors in collaboration with 
M. Nakamaye [BGN1,BGN2] were then able to construct many \Se metrics on certain 
connected sums of $S^2\times S^3$ as well as modify the Johnson-Koll\'ar computer code 
to handle more general log del Pezzo surfaces with higher Fano index. The original 
Johnson-Koll\'ar list contained several examples where the existence of a K\"ahler-Einstein 
metric was still in question.  One of these was treated in [BGN2] while two more were 
handled recently by C. Araujo [Ar]. It is the purpose of this note to show that the two new 
anticanonically embedded log del Pezzo surfaces shown in [Ar] to admit K\"ahler-Einstein 
metrics can be used to construct families of new \Se metrics on $8\#(S^2\times S^3).$ 

It is well-known [FK, BFGK] 
that $8\#(S^2\times S^3)$ admits regular \Se metrics; in fact, it 
admits an 8 complex dimensional family of \Se metrics [BGN1]. But up until now it was not 
known whether there are any non-regular \Se metrics on $8\#(S^2\times S^3).$ In this 
note we prove the following:

\noindent{\sc Theorem A}: \tensl The 5-manifold $8\#(S^2\times S^3)$ admits two 
families of non-regular \Se metrics. Each family depends on 8 complex parameters. Hence, 
$8\#(S^2\times S^3)$ admits 3 distinct 8 complex parameter families of \Se 
metrics, one regular and two non-regular families.  
These metrics are all inequivalent as Riemannian metrics. 
\tenrm

\bigskip
\bigskip
\baselineskip = 10 truept
\centerline{\bf 1. Sasakian Structures on Links of Isolated Hypersurface
Singularities}
\bigskip
Although the purpose of this note is to describe the \Se geometry of 
$8\#(S^2\times S^3)$ we shall begin with
a very brief summary of the Sasakian and \Se geometry of links of isolated
hypersurface singularities defined by weighted homogeneous polynomials.
For more details we refer the reader to [BG2, BGN1, BGN2].
Consider the affine space $\bbc^{n+1}$ together with a weighted
$\bbc^*$-action $\bbc^*_\bfw$ given by $(z_0,\ldots,z_n)\mapsto
(\grl^{w_0}z_0,\ldots,\grl^{w_n}z_n),$ where the {\it weights} $w_j$ are
positive integers. It is convenient to view the weights as the components of a
vector $\bfw\in (\bbz^+)^{n+1},$ and we shall assume that
$\gcd(w_0,\ldots,w_n)=1.$ Let $f$ be a quasi-homogeneous polynomial, that is
$f\in \bbc[z_0,\ldots,z_n]$ and satisfies
$$f(\grl^{w_0}z_0,\ldots,\grl^{w_n}z_n)=\grl^df(z_0,\ldots,z_n),
\leqno{1.1}$$
where $d\in \bbz^+$ is the degree of $f.$ We are interested in the {\it
weighted affine cone} $C_f$ defined by
the equation $f(z_0,\ldots,z_n)=0.$ We shall assume that the origin in
$\bbc^{n+1}$ is an isolated singularity, in fact the only singularity, of
$f.$ Then the link $L_f$ defined by
$$L_f= C_f\cap S^{2n+1}, \leqno{1.2}$$
where
$$S^{2n+1}=\{(z_0,\ldots,z_n)\in \bbc^{n+1}|\sum_{j=0}^n|z_j|^2=1\}$$
is the unit sphere in $\bbc^{n+1},$ is a smooth manifold of dimension $2n-1.$
Furthermore, it is well-known [Mil] that the link $L_f$ is $(n-2)$-connected.

Recall that a Sasakian structure consists of a quadruple $\cals =(\xi,\eta,\Phi,g)$ where $g$ is 
a Riemannian metric, $\xi$ is a unit length Killing vector field, $\eta$ is a contact 1-form such 
that $\xi$ is its Reeb vector field, and $\Phi$ is a $(1,1)$ tensor field which annihilates $\xi$ 
and describes an integrable complex structure on the contact vector bundle $\cald 
=\hbox{ker}~\eta.$ 
On $S^{2n+1}$ there is a well-known  ``weighted'' Sasakian structure
$\cals_\bfw=(\xi_\bfw,\eta_\bfw,\Phi_\bfw,g_\bfw)$, where the vector field
$\xi_\bfw$ is the infinitesimal generator of the circle subgroup $S^1_\bfw
\subset \bbc^*_\bfw.$ This Sasakian structure on $S^{2n+1}$ induces a
Sasakian structure, also denoted by
$\cals_\bfw$, on the link $L_f.$ (See
[YK, BG2, BGN1] for details.
The quotient space $\calz_f$ of $S^{2n+1}$ by  $S^1_\bfw,$ or
equivalently the space of leaves of the characteristic foliation $\calf_\xi$
of $\cals_\bfw,$  is a compact K\"ahler orbifold which is a projective
algebraic variety embedded in the weighted projective $\bbp(\bfw)=
\bbp(w_0,w_1,\cdots,w_n),$ in such a way that there is a commutative diagram
$$\matrix{L_f &\ra{2.5}& S^{2n+1}_\bfw&\cr
  \decdnar{\pi}&&\decdnar{} &\cr
   \calz_f &\ra{2.5} &\bbp(\bfw),&\cr}$$
where the horizontal arrows are Sasakian and K\"ahlerian embeddings,
respectively, and the vertical arrows are principal $S^1$ V-bundles and
orbifold Riemannian submersions.  

We are interested in deformations of the Sasakian structure that leaves the Reeb vector field 
invariant.  Such deformations of a
given Sasakian structure $\cals=(\xi,\eta,\Phi,g)$ are obtained by adding to
$\eta$ a  continuous one parameter family of basic1-forms $\grz_t.$
We require that the 1-form $\eta_t=\eta +\grz_t$  satisfy the conditions
$$\eta_0=\eta, \qquad \grz_0=0,\qquad \eta_t\wedge (d\eta_t)^n\neq 0~~
\forall~~ t\in [0,1]. \leqno{1.5}$$   
Since $\grz_t$ is basic $\xi$ is the Reeb (characteristic) vector field
associated to $\eta_t$ for all $t.$ This gives rise to a Sasakian structure
$\cals_t=(\xi,\eta_t,\Phi_t,g_t)$ for all $t\in
[0,1]$ that has the same underlying contact structure and the same
characteristic foliation. In general these structures are inequivalent and
the moduli space of Sasakian structures having the same characteristic vector
field is infinite dimensional.

Suppose now we have  a link $L_f$ with a given
Sasakian structure $(\xi,\eta,\Phi,g)$. When can we find a 1-form $\grz$ such
that the deformed structure $(\xi,\eta+\grz,\Phi',g')$ is Sasakian-Einstein? This is a
Sasakian version [BGN3] of the Calabi problem and its solution is equivalent to
solving the corresponding Calabi problem on the space of leaves $\calz_f$.
Since a \Se manifold necessarily has positive Ricci tensor, its
Sasakian structure is necessarily positive. This also implies that the
K\"ahler structure on $\calz_f$ be positive, i.e. $c_1(\calz_f)$ can be
represented by a positive definite $(1,1)$ form. In this case there are
well-known obstructions to solving to solving the Calabi problem.  These
obstructions for finding a solution to the Monge-Ampere equations involve
the non-triviality of certain {\it multiplier ideal sheaves} [Na, DK] associated
with effective anti-canonical $\bbq$-divisors on the space of leaves $\calz_f.$
Consequently, if one can show that these multiplier ideal sheaves coincide
with the full structure sheaf, one obtains the existence of a positive
K\"ahler-Einstein metric on $\calz_f$ and hence, a \Se metric on $L_f.$

\bigskip
\baselineskip = 10 truept
\centerline{\bf 2. The Construction}
\bigskip

The proof of Theorem A is based first on previous work of the authors (and M. Nakamaye) 
[BG1,BGN1], which in turn is based on the work of Demailly and Koll\'ar [DK] and Johnson 
and Koll\'ar [JK], and second on the recent 
work of C. Araujo [Ar].  Indeed, in [Ar] it is shown 
that the two log del Pezzo surfaces $\calz_{10}\subset \bbp(1,2,3,5)$ and $\calz_{15}\subset 
\bbp(1,3,5,7)$ admit K\"ahler-Einstein metrics. The weighted homogeneous polynomials 
$f_{10}$ and $f_{15}$ describing these log del Pezzo surfaces consist of 17 and 19 
monomials, respectively. It is assumed that the coefficients of the monomials are such that 
the del Pezzo surface is quasi-smooth. Moreover, the coefficient of the $z_1z_2z_3$ term in 
the $\calz_{15}$ surface is non-vanishing, 
for otherwise it is not known whether $\calz_{15}$ 
admits a K\"ahler-Einstein metric [Ar]. These log del Pezzo surfaces are orbifolds and the 
total space of the circle V-bundles over them are the links of isolated hypersurface 
singularities in $\bbc^4$ described by the corresponding weighted homogeneous polynomial. 
Combining our previous work [BG1,BG2] with Theorem 4.1 of [Ar] gives

\noindent{\sc Lemma 2.1}: \tensl The total spaces $L_{10}$ and $L_{15}$ of the V-bundles 
over $\calz_{10}$ and $\calz_{15},$ respectively  whose first Chern classes are represented 
by the K\"ahler form on the corresponding log del Pezzo surface admit  \Se metrics. 
\tenrm

Next we need to identify the links $L_{10}$ and  $L_{15}.$ First we mention that $L_{10}$ 
and  $L_{15}$ are inequivalent as links. This follows by computing their Milnor numbers 
$$\mu =\mu(L_d)=\prod_{i=0}^4\bigl({d\over w_i}-1\bigr).$$ 
We find $\mu(L_{10})=84$ and $\mu(L_{15})=128.$

\noindent{\sc Lemma 2.2}: \tensl The links $L_{10}$ and  $L_{15}$ are diffeomorphic to 
$8\#(S^2\times S^3).$ \tenrm

\noindent{\sc Proof}: The result was stated without proof in [BGN1]. First we notice that in 
both cases the weights are pairwise relatively prime, so the hypersurfaces $\calz_{10}$ and 
$\calz_{15}$ are well-formed. Thus, by Lemma 5.8 of [BG2] there is no torsion in 
$H_2(L_{15},\bbz).$ Moreover, it is well known that the links of an isolated hypersurface are 
always simply connected in this dimension. Then since by Lemma 1 $L_{10}$ and
$L_{15}$ admit a \Se metric, it follows from Corollary 2.1.6 of [BG1] they are both spin. It then 
follows from a well known theorem of Smale [Sm] that both $L_{10}$ and $L_{15}$ must be 
of the form $k\#(S^2\times S^3)$ for some positive integer $k.$ To find $k$ we compute the 
second Betti number by the method of Milnor and Orlik [MO]. We give the details for the case 
$L_{15}$ as $L_{10}$ is similar. Let $\grD(t)$ denote the Alexander 
polynomial of the link $L_{15},$ and let $\Div\grD(t)$ denote its divisor in the group ring 
$\bbz[\bbc^*].$ Then $\Div\grD(t)$ takes the form 
$$\Div\grD(t) = 1+\sum_ia_i\grL_i$$
for some $a_i\in \bbz$ where $\grL_i=\Div (t^i-1).$ Then the second Betti number is given 
by
$$b_2(L_f)= 1+\sum_ia_i.$$
The Milnor-Orlik procedure gives
$$\eqalign{\Div\grD(t)=&(\grL_{15}-1)({\grL_{15}\over 7}-1)(\grL_5-1)(\grL_3-1)\cr
                                   =&(\grL_{15}+1)(\grL_{15}-\grL_5-\grL_3+1)\cr
                                    =&9\grL_{15}-\grL_5-\grL_3+1.}$$
Thus, $k=b_2(L_{15})=8.$ Here we have used the relations 
$\grL_a\grL_b=\gcd(a,b)\grL_{lcm(a,b)}.$ \hfill\za

\bigskip
\baselineskip = 10 truept
\centerline{\bf 3. Moduli of Sasakian-Einstein Structures on $8\#(S^2\times S^3)$}
\bigskip

Next we determine the effective parameters in the families of \Se metrics.
Recall that
$\bbp(\bfw)$ can be
defined as a scheme ${\rm Proj}(S(\bfw))$,
where
$$S(\bfw)=\bigoplus_dS^d(\bfw)=\bbc[z_0,z_1,z_2,z_3].$$
The ring of polynomials $\bbc[z_0,z_1,z_2,z_3]$ is graded with grading
defined by the weights $\bfw=(w_1,w_1,w_2,w_3)$. As a projective variety
we can embed $\bbp(\bfw)\subset \bbc\bbp^N$ and then the
group  $\gG_\bfw$ is a subgroup of $PGL(N,\bbc)$. Note that $\bbp(\bfw)$ is
a toric variety and we can describe its complex
automorphism group $\gG_\bfw$ explicitly as follows:
Let $\bfw=(w_0,w_1,w_2,w_3)$ be ordered with $w_0\leq w_1\leq w_2\leq w_3$.
We consider the group $G(\bfw)$ of automorphisms of the graded ring $S(\bfw)$
defined on generators by
$$\varphi_\bfw\pmatrix{z_0\cr z_1\cr z_2\cr z_3\cr}=
\pmatrix{f_0^{(w_0)}(z_0,z_1,z_2,z_3)\cr
f_1^{(w_1)}(z_0,z_1,z_2,z_3)\cr f_2^{(w_2)}(z_0,z_1,z_2,z_3)\cr
f_3^{(w_3)}(z_0,z_1,z_2,z_3)\cr},$$
where $f_i^{(w_i)}(z_0,z_1,z_2,z_3)$ is an arbitrary weighted homogeneous
polynomial of degree $w_i$ in $(z_0,z_1,z_2,z_3)$. This is
a finite dimensional Lie group and
it is a subgroup of $GL(N,\bbc)$.
Projectivising, we get $\gG_\bfw=\bbp_\bbc(G(\bfw))$.

Note that when $\bfw=(1,1,1,1)$ then $G(\bfw)=GL(4,\bbc)$.  Other than this
case three weights are never the same if $\bbp(\bfw)$ is well-formed. If two
weights coincide then $G(\bfw)$ contains $GL(2,\bbc)$ as a subgroup. Finally,
when all weights are distinct we can write
$$\varphi_\bfw\pmatrix{z_0\cr z_1\cr z_2\cr z_3\cr}=
\pmatrix{a_0z_0\cr
a_1z_1+f_1^{(w_1)}(z_0)\cr a_2z_2+f_2^{(w_2)}(z_0,z_1)\cr
a_3z_3+f_3^{(w_3)}(z_0,z_1,z_2)\cr}$$
where $(a_0,a_1,a_2,a_3)\in(\bbc^*)^4$ and $f_i^{(w_i)},\ \ i=1,2,3$
are weighted homogeneous polynomials of degree $w_i$. The simplest situation
occurs when $f_1=f_2=f_3$ are forced to vanish. Then $\gG_\bfw=(\bbc^*)^3$ is
the smallest it can possibly be as $\bbp(\bfw)$ is toric.

Let $S_\bfw^d\subset S(\bfw)$ be the vector subspace spanned by all monomials
in $(z_0,z_1,z_2,z_3)$ of degree $d=|w|-I$, and
let $\hat{S}^d(\bfw)\subset S^d(\bfw)$ denote subset all quasi-smooth
elements.
Then we define $m_\bfw^d$ to be the dimension of the subspace
generated by $\hat{S}_\bfw^d.$ Now the automorphism group $G(\bfw)$ acts on
$S_\bfw^d$ leaving the subset $\hat{S}^d(\bfw)$ of quasi-smooth polynomials
invariant. Thus, for each log del Pezzo surface
we define the moduli space
$$\calm_\bfw^d=\hat{S}_\bfw^d/G(\bfw)=\bbp(\hat{S}_\bfw^d)/\gG_\bfw,
$$
with $n_\bfw^d={\rm dim}(\calm_\bfw^d).$ Now there is an injective map
$$\calm_\bfw^d\ra{2.0} \calm^\bbc(\calz_\bfw),$$
and each element in $\calm_\bfw^d$
corresponds to a unique homothety class of K\"ahler-Einstein metrics modulo
$\gG_\bfw$ and hence, to a unique \Se structure on the corresponding 5-manifold
modulo the group $\gG_\bfw$ acting as CR automorphisms. 

Let us now consider our $\calz_{10}$. It is easy to see that
$S^{10}(1,2,3,5)$ is isomorphic to $\bbc^{20}$ 
and it is spanned by the monomials
$z_3^2$, $z_1z_2z_3$, $z_1^2z_2^2$, $z_1^5$, $z_0z_2^3$,
$z_0z_1^2z_3$, $z_0z_1^3z_2$, $z_0^2z_2z_3$, $z_0^2z_1z_2^2$, $z_0^2z_1^4$,
$z_0^3z_1z_3$, $z_0^3z_1^2z_2$, $z_0^4z_2^2$, $z_0^4z_1^3$, $z_0^5z_3$,
$z_0^{5}z_1z_2$, $z_0^{6}z_1^2$, $z_0^{7}z_2$, $z_0^{8}z_1,$ $z_0^{10}$. We take the
open submanifold $\hat S^{10}(1,2,3,5)\subset S^{10}(1,2,3,5)$.
This is acted on by the complex automorphism group, namely the group
$G(1,2,3,5)$ generated by
$$\varphi_\bfw\pmatrix{z_0\cr z_1\cr z_2\cr z_3\cr}=
\pmatrix{a_0^1z_0\cr
a_1^1z_1+a_1^2z_0^2\cr a_2^1z_2+a_2^2z_0^2+a_2^3z_0z_1\cr
a_3^1z_3+a_3^2z_0^5+a_3^3z_0^3z_1
+a_3^3z_0^2z_2+a_3^4z_0z_1^2+a_3^5z_2z_1\cr},$$
where $a_i^1\in\bbc^*$ and all other coefficients are in $\bbc.$
$G(1,3,5,8)$ is a 
12-dimensional complex Lie group acting on the open submanifold
$\hat S^{10}(1,3,5,8)\subset S^{10}(1,3,5,8)\approx\bbc^{20}.$ It follows 
that the quotient is an 8-dimensional complex manifold which by the
Bando-Mabuchi Theorem [BM] is the moduli space of positive K\"ahler-Einstein
metrics on the underlying compact orbifold $\calz_{10}$.

Similarly for $\calz_{15}$ it is easy to see that
$S^{15}(1,3,5,7)$ is isomorphic to $\bbc^{19}$
and that it is spanned by the monomials
$z_2^3$, $z_1z_2z_3$, $z_1^5$, $z_0z_3^2$,
$z_0z_1^3z_2$, $z_0^2z_1z_2^2$, $z_0^2z_1^2z_3$, $z_0^3z_2z_3$, $z_0^3z_1^4$,
$z_0^4z_1^2z_2$, $z_0^5z_2^2$, $z_0^5z_1z_3$, $z_0^6z_1^3$, $z_0^7z_1z_2$,
$z_0^{8}z_3$, $z_0^{9}z_1^2$, $z_0^{10}z_2$, $z_0^{12}z_1$, $z_0^{15}.$ The
open submanifold $\hat S^{15}(1,3,5,7)\subset S^{15}(1,3,5,7)$ of quasi-smooth elements
is under the action of the complex automorphism group, namely the group
$G(1,3,5,7)$ generated by
$$\varphi_\bfw\pmatrix{z_0\cr z_1\cr z_2\cr z_3\cr}=
\pmatrix{a_0^1z_0\cr
a_1^1z_1+a_1^2z_0^3\cr a_2^1z_2+a_2^2z_0^5+a_2^3z_0^2z_1\cr
a_3^1z_3+a_3^2z_0^7
+a_3^3z_1z_0^4+a_3^4z_2^2z_0^2+a_3^5z_1^2z_0\cr},$$
where $a_i^1\in\bbc^*$ and all other coefficients are in $\bbc.$
$G(1,3,5,7)$ is a
11-dimensional complex Lie group acting on the open submanifold
$\hat S^{15}(1,3,5,7)\subset S^{15}(1,3,5,7)\approx\bbc^{19}.$ It follows
that the quotient is an 8-dimensional complex manifold. However, because of the condition in 
[Ar] that the coefficient of the term $z_1z_2z_3$ in the polynomial $f_{15}$ be 
non-vanishing, 
it is not known whether all the elements of $\hat S^{15}(1,3,5,7)$ admit a K\"ahler-Einstein 
metric.  Nevertheless, it is shown in [Ar] that there is a positive K\"ahler-Einstein metric on the 
open submanifold ${\tilde S}^{15}(1,3,5,7)$ defined by demanding that the coefficient of 
the $z_1z_2z_3$ term is not zero.  It then follows from [BGN1] that the moduli space has 
complex dimension 8. 

Thus for both $f_{10}$ and $f_{15}$ we  obtain an 8-complex dimensional family of \Se 
metrics on $8\#(S^2\times S^3).$ So the moduli space of \Se metrics on $8\#(S^2\times 
S^3)$ has at least three 8-complex dimensional families of \Se metrics one of which is 
regular and the 
other two non-regular. In a forthcoming work we show that these three families belong to 
distinct components. It also follows as in [BGN1] that the metrics are inequivalent as 
Riemannian metrics. This completes the proof of Theorem A. \hfill\za

\bigskip
\medskip
\centerline{\bf Bibliography}
\medskip
\font\ninesl=cmsl9
\font\bsc=cmcsc10 at 10truept
\parskip=1.5truept
\baselineskip=11truept
\ninerm

\item{[Ar]} {\bsc C. Araujo}, {\ninesl K\"ahler-Einstein Metrics for some quasi smooth 
log del Pezzo Surfaces}, preprint arXiv.math.AG/0111164.

\item{[BFGK]} {\bsc H. Baum, T. Friedrich, R. Grunewald, I. Kath}, {\ninesl
Twistors and Killing spinors on Riemannian manifolds},
Teubner-Texte zur Mathematik, 124. 
B. G. Teubner Verlagsgesellschaft mbH, Stuttgart, 1991.

\item{[BG1]} {\bsc C. P. Boyer and  K. Galicki}, {\ninesl On Sasakian-Einstein
Geometry}, Int. J. Math. 11 (2000), 873-909.

\item{[BG2]} {\bsc C. P. Boyer and  K. Galicki}, {\ninesl New Einstein Metrics
in Dimension Five}, J. Diff. Geom. 57 (2001) 443-463.
 
\item{[BGN1]} {\bsc C. P. Boyer, K. Galicki, and M. Nakamaye}, {\ninesl On the
Geometry of Sasakian-Einstein 5-Manifolds}, math.DG/0012047; to appear in Math. Ann.

\item{[BGN2]} {\bsc C. P. Boyer, K. Galicki, and M. Nakamaye}, {\ninesl
Sasakian-Einstein Structures on $\scriptstyle{9\#(S^2\times S^3)}$}, Trans. 
Amer. Math. Soc. 354 (2002), 2983-2996; math.DG/0102181.

\item{[BGN3]} {\bsc C. P. Boyer, K. Galicki, and M. Nakamaye}, {\ninesl On 
Positive Sasakian Geometry}, to appear in Geometriae Dedicata; math.DG/0104126.

\item{[DK]} {\bsc J.-P. Demailly and J. Koll\'ar}, {\ninesl Semi-continuity of
complex singularity exponents and K\"ahler-Einstein metrics on Fano
orbifolds}, preprint AG/9910118, Ann. Scient. Ec. Norm. Sup. Paris 34 (2001), 525-556.

\item{[FK]} {\bsc T. Friedrich and I. Kath}, {\ninesl Einstein manifolds of
dimension five with small first eigenvalue of the Dirac operator}, J. Diff.
Geom. 29 (1989), 263-279.

\item{[JK1]} {\bsc J.M. Johnson and J. Koll\'ar}, {\ninesl K\"ahler-Einstein
metrics on log del Pezzo surfaces in weighted projective 3-space}, Ann. Inst. 
Fourier 51(1) (2001) 69-79.
\item{[Mil]} {\bsc J. Milnor}, {\ninesl Singular Points of Complex
Hypersurfaces}, Ann. of Math. Stud. 61, Princeton Univ. Press, 1968.
\item{[Na]} {\bsc A.M. Nadel}, {\ninesl Multiplier ideal sheaves
and existence of K\"ahler-Einstein metrics of positive scalar curvature}, Ann.
Math. 138 (1990), 549-596.
\item{[MO]} {\bsc J. Milnor and P. Orlik}, {\ninesl Isolated singularities
defined by weighted homogeneous polynomials}, Topology 9 (1970), 385-393.

\item{[Sm]} {\bsc S. Smale}, {\ninesl On the structure of 5-manifolds},
Ann. Math. 75 (1962), 38-46.
\item{[YK]} {\bsc K. Yano and M. Kon}, {\ninesl
Structures on manifolds}, Series in Pure Mathematics 3,
World Scientific Pub. Co., Singapore, 1984.
\medskip
\bigskip \line{ Department of Mathematics and Statistics
\hfil August 2002} \line{ University of New Mexico \hfil }
\line{ Albuquerque, NM 87131 \hfil } \line{ email: cboyer@math.unm.edu,
galicki@math.unm.edu \hfill} \line{ web pages:
http://www.math.unm.edu/$\tilde{\phantom{o}}$cboyer,
http://www.math.unm.edu/$\tilde{\phantom{o}}$galicki \hfil}

\bye